\documentclass{preprint}

\ifoot{Submitted, July 2019} \ofoot{~} \ohead{\thepage}

\titlehead{Submitted July 2019}

\usepackage{amssymb,amsmath,graphicx,dsfont,xcolor,booktabs,hyperref,listings,pgf,subcaption,units,placeins}

\newcommand{\FL}{{\cal F}} \newcommand{\SO}{{\cal S}}
\newcommand{\IN}{{\cal I}}

\newcommand{\Ft}{\mathbf{F}} \newcommand{\Et}{\mathbf{E}}

\newcommand{\sigmat}{\boldsymbol{\sigma}}
\newcommand{\Sigmat}{\boldsymbol{\Sigma}}

\renewcommand{\div}{\operatorname{div}}
\newcommand{\tr}{\operatorname{tr}}
\newcommand{\id}{\operatorname{id}} \newcommand{\Vt}{{\mathbf V}}
\newcommand{\Lt}{{\mathbf L}} \newcommand{\Wt}{{\mathbf W}}

\newcommand{\vt}{\mathbf{v}} 
 
 \newcommand{\ft}{\mathbf{f}}
 
\newcommand{\ut}{\mathbf{u}} 
\newcommand{\nt}{\vec{n}}

\newcommand{\PC}{C} \newcommand{\pf}{F}

\title{Parallel time-stepping for fluid-structure interactions}
\author{Nils Margenberg \and Thomas Richter\thanks{Otto-von-Guericke
    Universit\"at Magdeburg, \url{nils.margenberg@ovgu.de},
    \url{thomas.richter@ovgu.de}}}

\begin{document}

\maketitle

\begin{abstract}
  We present a parallel time-stepping method for fluid-structure
  interactions. The interaction between the incompressible
  Navier-Stokes equations and a hyperelastic solid is formulated in a
  fully monolithic framework. Discretization in space is based on
  equal order finite element for all variables and a variant of the
  Crank-Nicolson scheme is used as second order time integrator. To
  accelerate the solution of the systems, we analyze a parallel-in
  time method. For different numerical test cases in 2d and in 3d we
  present the efficiency of the resulting solution approach. We also
  discuss some special challenges and limitations that are connected
  to the special structure of fluid-structure interaction problem.
\end{abstract}

\section{Introduction}

Fluid structure interactions appear in various problems ranging from
classical applications in engineering like the design of ships or
aircrafts, the design of wind turbines, but they are also present in
bio/medical systems describing the blood flow in the heart or in
general problems involving the cardiovascular system. The typical
challenge of fluid-structure interactions is two-fold. First, the
special coupling character that stems from the coupling of a
hyperbolic-type equation - the solid problem - with a parabolic-type
equation - the Navier-Stokes equations. Second, the moving domain
character brings along severe nonlinearities that have a non-local
character, as geometrical changes close to the moving fluid-solid
interface might have big impact on the overall solution.

Numerical approaches can usually be classified into \emph{monolithic
  approaches}, where the coupled fluid-structure interaction system is
taken as one entity and into \emph{partitioned approaches}, where two
separate problems - for fluid and solid - are formulated and where the
coupling between them is incorporated in terms of an outer (iterative)
algorithm. This second approach has the advantage that difficulties
are isolated and that perfectly suited numerical schemes can be used
for each of the subproblems. There are however application classes
where partitioned approaches either fail or lack efficiency. The
\emph{added mass effect}~\cite{CausinGerbeauNobile2005} exactly
describes this special stiffness connected to fluid-structure
interactions. It is typical for problems with similar densities in the
fluid and the solid - as it happens in the interaction of blood and
tissue or in the interaction of water and the solid structure of a
vessel. Here, monolithic approaches are considered to be favourable.

Monolithic approaches all give rise to strongly coupled, usually very
large and nonlinear algebraic systems of equations. Although there has
been substantial progress in designing efficient numerical schemes for
tackling the nonlinear
problems~\cite{HronTurek2006a,HeilHazelBoyle2008,FernandezGerbeau2009}
(usually by Newton's method) and the resulting linear
systems~\cite{GeeKuettlerWall2010,TurekHronMadlikRazzaqWobkerAcker2010,Richter2015,LangerYang2017,AulisaBnaBornia2018},
the computational effort is still immense. Numerically accurate
results and efficient approaches for 3d problems are still very
rare~\cite{Richter2017,FailerRichter2019}.

In this contribution we exploit the perspectives (and limitations) of
parallel time-stepping schemes for fluid-structure interaction
systems. We will have to face various difficulties that come from the
special type of equations and coupling, such as the hyperbolic
property of the solid problem and the saddle-point structure of the
fluid system.
Parallel time-stepping methods are well established in
the literature~\cite{fischer2005parareal, samaddar2010parallelization,
  croce2014parallel, haut2014asymptotic, blouza2011parallel,
  baudron2014parareal, kreienbuehl2015numerical},
but up to now there is little experience with
fluid-structure interactions.
However Parareal has shown to be unstable for hyperbolic problems,
where~\cite{ruprecht2018wave} found the phase error between the
fine and coarse propagators to be responsible, which already
foreshadows possible complications for our use case.

In the following section we will shortly describe a monolithic
Arbitrary Lagrangian Eulerian formulation for fluid-structure
interaction problems. Further, we detail on the discretization of this
system in space (with continuous finite elements) and time (with
classical time-stepping methods). The following 3rd section will focus
on describing a parallel time-stepping scheme and the special
requirements for a realization in terms of fluid-structure
interactions. Finally, in section~\ref{sec:num} we present numerical
results and discuss possible applications but we also address some
limitations.

\section{Fluid-structure interactions}\label{sec:fsi}

The presentation within this section mainly follows~\cite{Richter2017}.
We consider fluid-structure interaction problems coupling an
incompressible fluid with a hyperelastic solid. By $\hat\SO$ we denote
the Lagrangian reference framework of the solid, by $\SO(t)$ its
current configuration in Eulerian coordinates. By $\FL(t)$ we denote
the fluid domain at time $t$ matching the solid at the common
interface $\IN(t)=\partial\SO(t)\cap \partial\FL(t)$. By
$\Omega(t):=\FL(t)\cup\IN(t)\cup \SO(t)$ we denote the Eulerian
fluid-structure interaction domain. The domains $\FL(t),\SO(t)$ and
$\Omega(t)$ are all either two-dimensional or three-dimensional. The
boundary of the fluid domain
$\partial\FL(t)=\IN(t)\cup \Gamma_f^D(t)\cup\Gamma_f^{out}(t)$ is split
into the interface, a Dirichlet part $\Gamma_f^D(t)$ (usually inflow
or rigid walls) and an outflow part $\Gamma_f^{out}(t)$, where we ask
for the do-nothing condition~\cite{HeywoodRannacherTurek1992}. For
simplicity we assume Dirichlet conditions at the solid boundary apart
from the interface $\partial\SO(t)=\IN(t)\cup\Gamma_s^D(t)$. Finally,
by $I=[0,T]$ we denote the time interval. With the density $\rho_f$
the incompressible Navier-Stokes equations are given by
\begin{equation}\label{fluid}
  \begin{aligned}
    \div\,\vt_f=0,\quad \rho_f\big(\partial_t\vt_f +
    (\vt_f\cdot\nabla)\vt_f\big) -\div\,\sigmat_f(\vt_f,p_f)
    &=\rho_f\ft_f,
    &\text{in }&I\times \FL(t)\\
    \vt_f&=\vt_f^D&\text{on }&I\times \Gamma_f^D(t)\\
    \rho_f\nu_f\nabla\vt_f\nt_f - p_f\nt_f&=0&\text{on }&I\times
    \Gamma_f^{out}(t)\\
    \vt_f&=d_t \hat \ut_s&\text{on
    }&I\times \IN(t)\\
    \vt_f&=\vt_f^0&\text{on }&\{0\}\times \FL(0),
  \end{aligned}
\end{equation}
with the right hand side field $\ft_f$, boundary data $\vt_f^D(t)$ and
the interface velocity $d_t\ut_s$ coming from the coupling to the
solid equation. By $\vt_f^0$ we denote the initial velocity. The solid
problem in term is given in Lagrangian reference formulation on
$\hat\SO$ as
\begin{equation}\label{solid}
  \begin{aligned}
    \hat\rho_sd_{tt}\hat\ut_s - \widehat{\div}\big(
    \hat\Ft_s\hat\Sigmat_s(\hat \ut_s)\big)&=\hat\rho_s \hat\ft_s
    &\text{in }&I\times \hat \SO\\
    \hat \ut_s&=\hat \ut_s^D&\text{on }&I\times \hat\Gamma_s^D\\
    \hat\Ft_s\hat\Sigmat_s(\hat\ut_s)\hat\nt_s &=
    \sigma_f(\vt_f,p_f)\nt_f&\text{on }&I\times \hat \IN,\\
    d_t \hat\ut_s=\hat\vt_s^0,\quad \ut_s&=\hat\ut_s^0 &\text{on }
    &\{0\}\times\hat\SO,
  \end{aligned}
\end{equation}
where we denote by $\hat\Ft_s = I+\hat\nabla\hat\ut_s$ the deformation
gradient, by $\hat\rho_s$ the reference density, by $\hat\ft_s$ the
right hand side vector, boundary data by $\hat\ut_s^D$ and by
$\sigmat_f(\vt_f,p_f)\nt_f$ the normal stresses from the coupling to
the fluid equations. The hats ``$\hat{\,\cdot\,}$'' are added to
distinguish Lagrangian variables form their counterparts in the
reference framework. The attribution of the kinematic interface
condition $\vt_f=d_t\ut_s$ to the fluid problem and the dynamic
condition to the solid problem is artificial, as the coupled system
of~(\ref{fluid}) and~(\ref{solid}) must be considered as one entity.
Finally, as material models we consider a Newtonian fluid and a St.
Venant Kirchhoff solid
\begin{equation}\label{stress}
  \sigmat_f(\vt_f,p_f) = \rho_f\nu_f(\nabla\vt_f+\nabla\vt_f^T)-p_f
  I,\quad
  \hat\Sigmat_s = 2\mu_s \hat\Et_s +\lambda_s\tr(\hat\Et_s)I,
\end{equation}
where we denote by $\hat\Et_s=\frac{1}{2}(\hat\Ft_s^T\hat\Ft_s-I)$ the
Green-Lagrange strain tensor and by $\mu_s,\lambda_s$ the Lam\'e
parameters, by $\nu_f$ the kinematic viscosity.
By~(\ref{fluid}),~(\ref{solid}) and~(\ref{stress}) we denote the
fluid-structure interaction problem.

\subsection{Arbitrary Lagrangian Eulerian coordinates}

To overcome the discrepancy between the Eulerian fluid framework and
the Lagrangian solid framework, we map the flow problem to a fixed
reference domain $\hat\FL$ that fits the Lagrangian solid domain
$\hat\SO$. Here we assume that $\hat\FL=\FL(0)$ is just the known
fluid domain at initial time $t=0$. By $\hat T_f(t):\hat\FL\to \FL(t)$
we denote the reference map, by $\hat\Ft_f:=\hat\nabla\hat T_f$ its
gradient and by $\hat J_f:=\det\,\Ft_f$ its determinant. By
$\hat\vt_f(\hat x,t) = \vt_f(x,t)$ with $x=\hat T_f(x,t)$ we denote
the ALE representation of the Eulerian variable $\vt_f$ (same for the
pressure). The Arbitrary Lagrangian Eulerian formulation (ALE) goes
back to the
70s~\cite{HirtAmsdenCook1974,HughesLiuZimmermann1981,Donea1982}, and
usually consists of adding convective terms with respect to the motion
of the domain. We use a strict mapping to the fixed reference system
and formulate the coupled variational formulation on this arbitrary
framework $\hat\FL$. Details are given
in~\cite{RichterWick2010,Richter2017}. To close the ALE-formulation we
construct the map $\hat T_f:=\id + \hat\ut_f$ by means of an extension
of the solid deformation to the fluid domain, denoted by $\hat \ut_f$.
For small deformations of the fluid domain a simple harmonic extension
is sufficient
\[
  -\hat\Delta \hat\ut_f =0\text{ in }\hat\FL,\quad
  \hat\ut_f=\hat\ut_s\text{ on }\hat\IN,\quad \hat\ut_f=0\text{ on
  }\partial\hat\FL\setminus \hat\IN,
\]
while  problems with large changes in the fluid domain require more care
in extending the deformation. Details are given in~\cite[Section
5.3.5]{Richter2017} and the references therein. We finally give the
complete system. From hereon, all problems are given in the reference
system such that we skip the hat for better readability
\begin{equation}\label{ALE}
  \begin{aligned}
    J_f\Ft_f^{-1}:\nabla \vt^T = 0,\quad \rho_fJ_f\big(\partial_t\vt_f
    +
    \nabla\vt_f\Ft_f^{-1}(\vt_f-\partial_t\ut_f)\big)\qquad\\
    -\div\big( J\sigmat_f(\vt_f,p_f)\Ft_f^{-T}\big) &= \rho_f\ft_f
    &   \text{in }&I\times \FL\\
    \partial_t \ut_s =\vt_s,\quad \rho_s \partial_t\vt_s - \div\big(
    \Ft_s\Sigmat_s(\ut_s)\big)&=\rho_s \ft_s
    &\text{in }&I\times  \SO\\
    \vt_f=\vt_s,\quad J_f\sigma_f(\vt_f,p_f)\Ft_f^{-T}\nt_f
    +\Ft_s\Sigmat_s(\ut_s)\nt_s&=0    &\text{on }&I\times \IN,\\
    \vt_f=\vt_f^D\text{ on }I\times \Gamma_f^D(t),\quad
    \rho_f\nu_fJ_f\nabla\vt_f\Ft_f^{-1} \Ft_f^{-T}\nt_f -
    J_f\Ft_f^{-T}p_f\nt_f&=0 &\text{on }&I\times
    \Gamma_f^{out}(t)\\
    \ut_s&=\ut_s^D&\text{on }&I\times \Gamma_s^D\\
    \vt_f=\vt_f^0\text{ on }\{0\}\times \FL(0),\quad
    \vt_s=\vt_s^0,\quad \ut_s&=\ut_s^0 &\text{on } &\{0\}\times\SO.
  \end{aligned}
\end{equation}
Note that $\nt_f$ and $\nt_s$ are the outward facing normals in the
reference framework such that the fluid stresses are given in terms of
the Piola transform. The fluid stress tensor in ALE formulation reads
\begin{equation}\label{ALE:stress}
  \sigmat_f(\vt_f,p_f) =
  \rho_f\nu_f(\nabla\vt_f\Ft_f^{-1}+\Ft_f^{-T}\nabla\vt_f^T)-p_f
  I,\quad
  \hat\Sigmat_s = 2\mu_s \Et_s +\lambda_s\tr(\Et_s)I,
\end{equation}
In~(\ref{ALE}) we have split the hyperbolic solid problem into a
system of first order (in time) equations by introducing the solid
velocity $\vt_s$ that - on the interface - matches the fluid velocity.

\subsection{Variational formulation and finite element discretization}

To prepare for a discretization with finite elements we briefly sketch
the variational formulation of the fluid-structure interaction
problem~(\ref{ALE}) that also embeds the coupling conditions in
a variational sense. Kinematic and geometric coupling conditions are
taken care of by choosing global function spaces for the velocity and
deformation, i.e.
\[
  \vt\in \vt^D + H^1_0(\Omega;\Omega^D)^d,\quad \ut\in \ut^D +
  H^1_0(\Omega;\Omega^D)^d,
\]
where $d\in\{2,3\}$ is the spatial dimension,
$\Omega=\FL\cup\IN\cup\SO$ and $\Omega^D = \Gamma_f^D\cup\Gamma_s^D$
the combined Dirichlet boundary and $\ut^D,\vt^D$ are extensions of
the Dirichlet data into the domain. Solid and fluid velocity (and
deformation) are defined as the restrictions of $\vt$ and $\ut$ to the
respective domain. The dynamic condition is realized by testing the
momentum equations for both subproblems by one common and continuous
(in the $H^1$-sense) test functions $\phi\in H^1_0(\Omega;\Omega^D)^d$
and summing up both equations.

\begin{equation}\label{variational}
  \begin{aligned}
    \big(J \Ft^{-1}:\nabla \vt^T,\xi\big)_\FL +
    \big(\rho_fJ\big(\partial_t\vt +
    \nabla\vt\Ft^{-1}(\vt-\partial_t\ut)\big),\phi\big)_{\FL}\quad \\
    +\big(J\sigmat_f(\vt,p)\Ft^{-T},\nabla\phi\big)_{\FL} -\big\langle
    \rho_f\nu_f J\Ft^{-T}\nabla \vt^T \Ft^{-T}\nt_f
    ,\phi \big\rangle_{\Gamma_f^{out}} \quad \\
    +\big(\rho_s \partial_t\vt,\phi\big)_{\SO} + \big(
    \Ft\Sigmat_s(\ut),\nabla\phi\big)_\SO
    &= \big(\rho_f\ft_f,\phi\big)_\FL + \big(\rho_s
    \ft_s,\phi\big)_\SO\\
    \big(\partial_t \ut -\vt,\psi_s\big)_{\SO}
    +\big(\nabla\ut,\nabla\psi_f\big)_\FL &=0,\\
  \end{aligned}
\end{equation}
for all
\begin{equation}\label{test}
  \phi\in H^1_0(\Omega;\Omega^D)^d,\quad
  \xi\in L^2(\FL),\quad
  \psi_s\in L^2(\SO)^d,\quad \psi_f\in H^1_0(\FL).
\end{equation}
By $\langle\cdot,\cdot\rangle_{\Gamma}$ we denote the $L^2$-inner
product on a boundary segment $\Gamma$. The boundary term is
introduced to correct the normal stresses to comply with the
do-nothing outflow condition~\cite{HeywoodRannacherTurek1992}.
Dirichlet boundary values are embedded in the trial spaces, the
initial conditions will be realized in the time stepping schemes.

Next, let $\Omega_h$ be a triangulation (or more general a mesh) of
the domain $\Omega$ satisfying:
\begin{itemize}
\item The elements $K\in\Omega_h$ are open polytopes (we consider
  quadrilaterals or hexahedras, but other element types are possible).
\item Two different elements $K,K'\in\Omega_h$ do not overlap, i.e.
  $K\cap K'=\emptyset$, their boundaries $\partial K\cap \partial K'$
  either don't overlap, or they either overlap in a common node, a
  complete common edge or (in 3d) a (complete) common face.
\item To allow for local mesh refinement we relax the previous
  assumption and allow that an edge (or face) of one element $K$ is
  met by the edges (or faces) of two (or four) refined elements.
\item We assume that all elements stem from the mapping of one
  reference element $\hat T_T:\hat K\to K$, where $\hat K$ is the
  unit-quad (or unit-cube, or unit-tetrahedra, ...) and we assume that
  \[
    \|\nabla^k\hat T_T\|_\infty\|\nabla^k\hat T_T^{-1}\|_\infty\le
    c,\quad k=0,1,
  \]
  is uniformly bounded in the mesh size $h_T=\operatorname{diam}(T)$
  to allow for standard interpolation estimates.
\item We assume that the interface $\IN$ in the ALE reference
  configuration in resolved by the mesh, i.e. for all elements
  $K\in\Omega_h$ it holds $K\cap \IN = \emptyset$.
\end{itemize}

These assumptions are slight variations of typical requests on the
structural regularity and the form regularity of finite element
meshes, see~\cite{Ciarlet1978}. Hanging nodes on 2:1 balanced meshes
are also well-established in literature~\cite{BeckerBraack2000a}.
Details on the specific implementation in the finite element library
\emph{Gascoigne 3D}\cite{Gascoigne3D} are given in~\cite[Section
4.2]{Richter2017}. The last assumption is important to guarantee good
approximation properties, as fluid-structure interactions are
\emph{interface problems}, where the solution has limited regularity
across the interface and non-fitted meshes give rise to a breakdown
in accuracy~\cite{FreiRichter2014,FreiRichter2017}.

The finite element discretization of~(\ref{variational}) is
straightforward. We choose continuous polynomial spaces assembled on
the mesh $\Omega_h$ as
\[
  V_h^{(r)} := \{\phi\in C(\bar\Omega)\;|\; \phi\big|_T\circ \hat
  T_T^{-1}\in Q^{(r)}\},
\]
where $Q^{(r)}$ is the space of bi- or tri-polynomial of degree $r$ on
quadrilateral or hexahedral meshes
\[
  Q^{(r)} = \{ x_1^{\alpha_1}\cdots x_d^{\alpha_d},\; 0\le
  \alpha_1,\dots,\alpha_d\le r\},
\]
and where the \emph{iso-parametric} reference map
$\hat T_T\in [Q^{(r)}]^d$ comes from this same space. The
\emph{iso-parametric} setup is able to give optimal order
approximations on domains with curved boundaries, see~\cite{Lenoir1986}
or \cite[Section 4.2.3]{Richter2017}. Throughout this paper we choose
discrete trial- and test-spaces for the discretization
of~(\ref{variational}),~(\ref{test}) as
\[
  \vt_h\in \vt^D_h + \Vt_h\quad (\Vt_h:=[V_h^{(2)}]^d),\quad \ut_h\in
  \ut^D_h + \Vt_h,\quad p_h\in Q_h:= V_h^{(2)},
\]
and
\[
  \Lt_h = [V_h^{(2)}]^d,\quad \Wt_h = [V_h^{(2)}]^d,
\]
with the necessary modifications to implement Dirichlet values.

As the equal order finite element pair $V_h\times Q_h$ does not fulfil
the inf-sup condition we we add additional stabilization terms of
local projection~\cite{BeckerBraack2001} or internal jump
type~\cite{BurmanHansbo2006}. Required small modifications in terms of
the ALE formulation in fluid-structure interactions are discussed
in~\cite{Molnar2015} or~\cite[Section 5.3.3]{Richter2017}. For
simplicity we assume that no mechanisms for stabilizing dominant
transport are required.

\subsection{Time discretization}

For temporal discretization in the time interval $I=[0,T]$, we start
by a splitting into discrete time-steps
\[
  0=t_0<t_1<\cdots < T_N = T,\quad k:=t_n-t_{n-1}.
\]
To simplify notation we assume that this distribution is uniform with
$k$ being the same in all time steps $t_{n-1}\mapsto t_n$.
Modifications are however straightforward. At time $t_n$ we denote by
$\vt_n=\vt_h(t_n)$, $\ut_n=\ut_h(t_n)$, $p_n=p_h(t_n)$ the discrete
approximations and by $\Ft_n:=I+\nabla\ut_n,\quad J_n:=\det\,\Ft_n$,
the deformation gradient and its determinant. By the index
$n-\frac{1}{2}$ we denote the mean values on $I_n=[t_{n-1},t_n]$, e.g.
$\vt_n:=(\vt_{n-1}+\vt_n)/2$ or $J_n:=(J_{n-1}+J_n)/2$.

As a compromise between accuracy (second order), very good stability
properties (globally A-stable), simplicity and efficiency (simple
one-step scheme) we consider an implicitly shifted version of the
Crank-Nicolson
method~\cite{LuskinRannacher1982,Rannacher1984},~\cite[Section
5.1.2]{Richter2017}.

Applied to the finite element approximation to~(\ref{variational}),
each time step $t_{n-1}\mapsto t_n$ of the fully discrete problem
reads:
\begin{multline}
  \big(\rho_fJ_{n-\frac{1}{2}} \big((\vt_n-\vt_{n-1}) -
  \nabla\vt_{n-\frac{1}{2}}\Ft_{n-\frac{1}{2}}^{-1}
  (\ut_n-\ut_{n-1})  \big),\phi\big)_{\FL} \quad \\
  +\big(\rho_s (\vt_n-\vt_{n-1}),\phi\big)_{\SO}
  +\big(\ut_n-\ut_{n-1},\psi_s\big)_{\SO}\quad \\
  +k A_{div}(\vt^n,\ut^n)(\xi) + k \theta
  A_{NS}(\vt^n,\ut^n,p^n)(\phi)
  + k (1-\theta) A_{NS}(\vt^{n-1},\ut^{n-1},p^n)(\phi) \quad \\
  + k \theta A_{ES}(\vt^n,\ut^n)(\phi)
  + k (1-\theta) A_{ES}(\vt^{n-1},\ut^{n-1})(\phi) \quad \\
  -k\theta\big(\vt_n,\psi_s\big)_{\SO}-k(1-\theta)\big(\vt_n,\psi_s\big)_{\SO}
  +k\theta\big(\nabla\ut_n,\nabla\psi_f\big)_\FL
  +k(1-\theta)\big(\nabla\ut_{n-1},\nabla\psi_f\big)_{\FL} = F(\phi)
\end{multline}
for all $\phi\in \Vt_h,\xi\in Q_h,\psi_s\in \Lt_h$ and
$\psi_f\in \Wt_h$, where for simplicity of notation we introduced the
notation
\begin{equation}\label{semiforms}
  \begin{aligned}
    A_{div}(\vt,\ut)(\xi) &=
    \big(J \Ft^{-1}:\nabla \vt^T,\xi\big)_\FL\\
    A_{NS}(\vt,\ut,p)(\phi) &= \big(\rho_fJ
    \nabla\vt\Ft^{-1}\vt,\phi\big)_{\FL}
    +\big(J\sigmat_f(\vt,p)\Ft^{-T},\nabla\phi\big)_{\FL} -\big\langle
    \rho_f\nu_f J\Ft^{-T}\nabla \vt^T \Ft^{-T}\nt_f
    ,\phi \big\rangle_{\Gamma_f^{out}} \quad \\
    A_{ES}(\vt,\ut)(\phi) &= \big(
    \Ft\Sigmat_s(\ut),\nabla\phi\big)_\SO\\
    F(\phi) &=
    k\big(\rho_f(\theta\ft_f(t_n)+(1-\theta)\ft_f(t_{n-1})),\phi\big)_\FL
    + k\big(\rho_s (\theta
    \ft_s(t_n)+(1-\theta)\ft_s(t_{n-1})),\phi\big)_\SO.
  \end{aligned}
\end{equation}
For $\theta=\frac{1}{2}$ this is the typical Crank-Nicolson scheme,
$\theta=1$ would give the backward Euler method of first order. We
usually consider $\theta=\frac{1}{2}+\theta_0k$ with a small parameter
$\theta_0$ and refer to Section~\ref{sec:num}. A proper combination of
three substeps with different choices of $\theta$ would give the
fractional step theta method which is of second order and strongly
A-stable, see~\cite{TurekRivkindHronGlowinski2006}. The discretization
of the nonlinear terms including time derivatives (e.g. coming from
the domain convection) is still discussed in literature. However,
different choices give similar stability and accuracy results,
see~\cite{RichterWick2015_time},~\cite[Section 5.1.2]{Richter2017}.

\section{Parallel time-stepping}\label{sec:par}
The motivation for increased parallelism in numerical algorithms arises
from the architecture of modern computers with an ever-increasing
number of processing units. Most numerical solvers for differential
equations are formulated as sequential algorithms, which denies
exploitation of parallelism.

The Parareal algorithm, which is probably the best known parallel time
stepping method, promises to bypass this problem by basically
providing a wrapper around common sequential algorithms. Although it
is known in its fundamentals since 2001, it hasn't been widely applied
to fluid structure interactions.

The Parareal-algorithm can be derived in multiple ways,
see~\cite{GanderVandewalle2007} for an overview. The first approach
introduced in~\cite{LionsMadayTurinici2001}, is based on a predictor
corrector scheme. First the time domain $I$ is divided into $L$
subintervals of equal length. Predictions with coarse timesteps $K$
are used to provide initial values for parallel computations with fine
timesteps $k\ll K$. These results are then used together with old
coarse predictions for correcting the solution, providing new initial
values. Timesteps are chosen such that
$\frac{T}{L}=m\cdot k=M\cdot K$. The algorithm is defined by using two
different propagators for the solution over one subinterval. While
$\PC$ is the propagator using the coarse time step size, $\pf$ is the
fine counterpart.

With these definitions and the remarks from the last chapter in mind,
the Parareal-algorithm for fluid-structure-interactions is defined by
a simple recursive formula
\begin{equation}\label{parareal:iteration}
  \begin{aligned}
    (\vt,\,\ut,\,p)^0_i &= (\vt^0,\,\ut^0)\\
    (\vt,\,\ut,\,p)^{l+1}_{i+1} &=
    \underbrace{\PC((\ut,\,\vt)^l_{i+1},\,t_{l+1},\,t_l)}_{\text{predictor}}
    + \underbrace{\pf((\ut,\,\vt)^l_i,\,t_{l+1},\,t_l) -
      \PC((\ut,\,\vt)^l_i,\,t_{l+1},\,t_l)}_{\text{corrector}}\,,
  \end{aligned}
\end{equation}
where $l\in\lbrace0,\dots,\,L-1\rbrace$ is the index of the
subinterval and $i\in\mathbb{N}$ is the iteration count. While the
predictor part is sequential (particularly the first iteration for
initialization), the fine propagations on each subinterval can be
parallelized.

\begin{figure}[h!]
  \begin{lstlisting}[basicstyle=\ttfamily,frame=tb,
    escapeinside={(*}{*)}]
    #pragma omp parallel
    {
      #pragma omp ordered nowait schedule(static)
      (*$(\ut,\,\vt)_{0}^\pf = (\ut,\,\vt)_{0}$*) for m = 0 to N-1
      #pragma omp ordered
      {
        omp_set_lock(&interval_locker[i])
        (*$(\ut,\,\vt)_{m}^\PC = \PC((\ut,\,\vt)_m,\,t_{m+1},\,t_m)$*)
        (*$(\ut,\,\vt)_{m}^{\textrm{final}} = (\ut,\,\vt)_{m+1}^\PC$*)
        (*$(\ut,\,\vt)_{m+1}^\pf = (\ut,\,\vt)_{m+1}^\PC$*)
        omp_unset_lock(&interval_locker[i]) }

      for k = 1 to K
        #pragma omp for nowait schedule(static)
        for m = 0 to N - k + 1
          omp_set_lock(&interval_locker[m])
          (*$(\ut,\,\vt)_{m+k}^\pf =
          \pf((\ut,\,\vt)_{m+1}^\pf,\,t_{m+k},\,t_{m+k-1})$*)
          omp_unset_lock(&interval_locker[m])

      (*$(\ut,\,\vt)_{0}^{\textrm{final}} = (\ut,\,\vt)_{0}^\pf$*)
      #pragma omp for ordered nowait schedule(static)
      for m = 0 to N -k
      #pragma omp ordered
      {
        omp_set_lock(&interval_locker[m])
        (*$(\ut,\,\vt)_{m+k-1}^\PC =
        (\ut,\,\vt)_{m+k-1}^{\textrm{final}}$*)
        omp_unset_lock(&interval_locker[m])
        (*$(\ut,\,\vt)_{m+k-1}^\PC =
        \PC((\ut,\,\vt)_{m+k-1}^\PC,\,t_{l+1},\,t_l)$*)
        omp_set_lock(&interval_locker[m + 1])
        (*$(\ut,\,\vt)_{m+k}^{\textrm{final}}=(\ut,\,\vt)_{m+k-1}^\PC
        + (\ut,\,\vt)_{m+k}^\pf - (\ut,\,\vt)_{m+k}^\PC$*)
        (*$(\ut,\,\vt)_{m+k}^\pf =
        (\ut,\,\vt)_{m+k}^{\textrm{final}}$*)
        omp_unset_lock(&interval_locker[m + 1])
      }
    }
  \end{lstlisting}
  \caption{Pseudo-code for Parareal algorithm with OpenMP directives,
    following the distributed task scheduling.}
  \label{fig:parareal-pseudocode}
\end{figure}

\FloatBarrier

Recently a weighted version of the Parareal scheme, the
$\theta$-Parareal scheme has been developed, introducing weights for
the contributions of the coarse propagators:
\begin{equation}\label{thetaparareal:iteration}
  \begin{aligned}
    (\vt,\,\ut,\,p)^0_i &= (\vt^0,\,\ut^0)\\
    (\vt,\,\ut,\,p)^{l+1}_{i+1} &=
    \theta^{l+1}_{i+1}\PC((\ut,\,\vt)^l_{i+1},\,t_{l+1},\,t_l) +
    \pf((\ut,\,\vt)^l_i,\,t_{l+1},\,t_l) -
    \theta^{l+1}_{i+1}\PC((\ut,\,\vt)^l_i,\,t_{l+1},\,t_l)
  \end{aligned}
\end{equation}
The idea here was to replace the predictor term by
\[\theta^{l+1}_{i+1}\PC((\ut,\,\vt)^l_{i+1},\,t_{l+1},\,t_l) +
  (1-\theta^{l+1}_{i+1})\PC((\ut,\,\vt)^l_i,\,t_{l+1},\,t_l)\,,\]
yielding the formula above. One idea for obtaining $\theta^{l}_{i}$ is
to minimize the discrepancy of fine and coarse solution:
\[\theta^{l}_{i+1}=
  \arg\limits_{\theta\in\mathbb{R}}\min \lVert
  \pf((\ut,\,\vt)^l_i,\,t_{l+1},\,t_l)
  -\PC((\ut,\,\vt)^l_{i+1},\,t_{l+1},\,t_l)\rVert\,.\] In which case
$\theta^{l}_{i}$ can be computed as
$\frac{\langle\pf,\,\PC\rangle}{\langle\PC,\,\PC\rangle}$ assuming the
same arguments as above. We will use
$\frac{\langle\pf,\,\PC\rangle}{\langle\PC,\,\PC\rangle
  \langle\pf,\,\pf\rangle}$, which is not a least squares solution,
but can be interpreted as angle penalization. The factors are computed
for each component separately and averaged to obtain a single scaling
factor. This way we make sure that $\theta^{l}_{i}\leq 1$ by
construction. Later on we will shortly investigate the different
approaches. In figure~\ref{fig:parareal-pseudocode} a shared memory
implementation with OpenMP is sketched, following the distributed task
scheduling from~\cite{Aubanel2011172}. This is nearly optimal and the
theoretically achievable speedup with $r=\frac{h_\pf}{h_\PC}$ is
\begin{equation}\label{eq:speedup}
  S = \frac{1}{r + \frac{K}{N}(1 + r)}.
\end{equation}
In our examples presented
next, we cannot rely on the fact
that $r$ is actually the ratio between the computational costs of the
coarse and fine propagator.
In section~\ref{sec:num} the speedup and the convergence of Parareal will be
investigated.

\section{Numerical examples}\label{sec:num}
In this section we will apply the Parareal algorithm to two FSI
problems. Further we discuss various issues that we were facing when
applying the Parareal algorithm to different configurations. Based on
these examples we will measure the error introduced by
the Parareal algorithm as well as the speedup and efficiency.
Measurements are taken on a two socket machine with two Intel Xeon
E5--2640 v4 each having 10 cores. The number of subintervals
in the Parareal algorithm are chosen as the total number of cores, 20.

\begin{figure}[t!]
  \begin{subfigure}[l]{\textwidth}
    \begin{center}
      \includegraphics[width=0.6\textwidth]{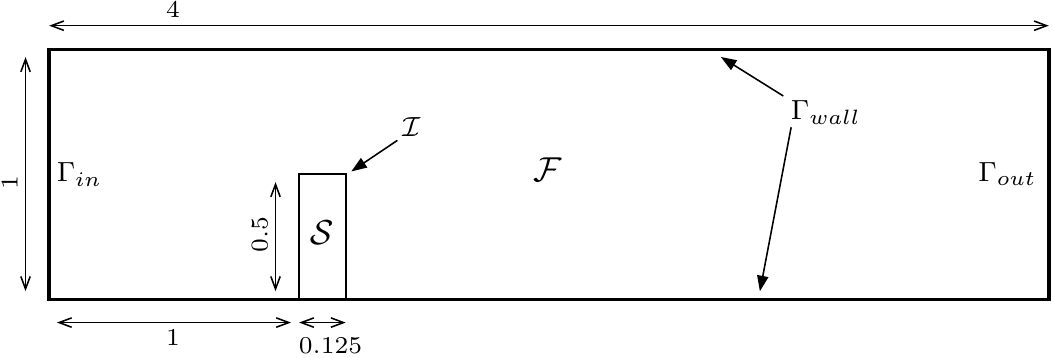}
    \end{center}
    \subcaption{2d configuration. Dirichlet inflow on $\Gamma_{in}$,
      do-nothing outflow condition on $\Gamma_{out}$ and no-slip
      Dirichlet condition on $\Gamma_{wall}$. By $\IN$ we denote the
      fluid-structure interface. }
  \end{subfigure}
  \begin{subfigure}[l]{\textwidth}
    \begin{center}
      \includegraphics[width=0.8\textwidth]{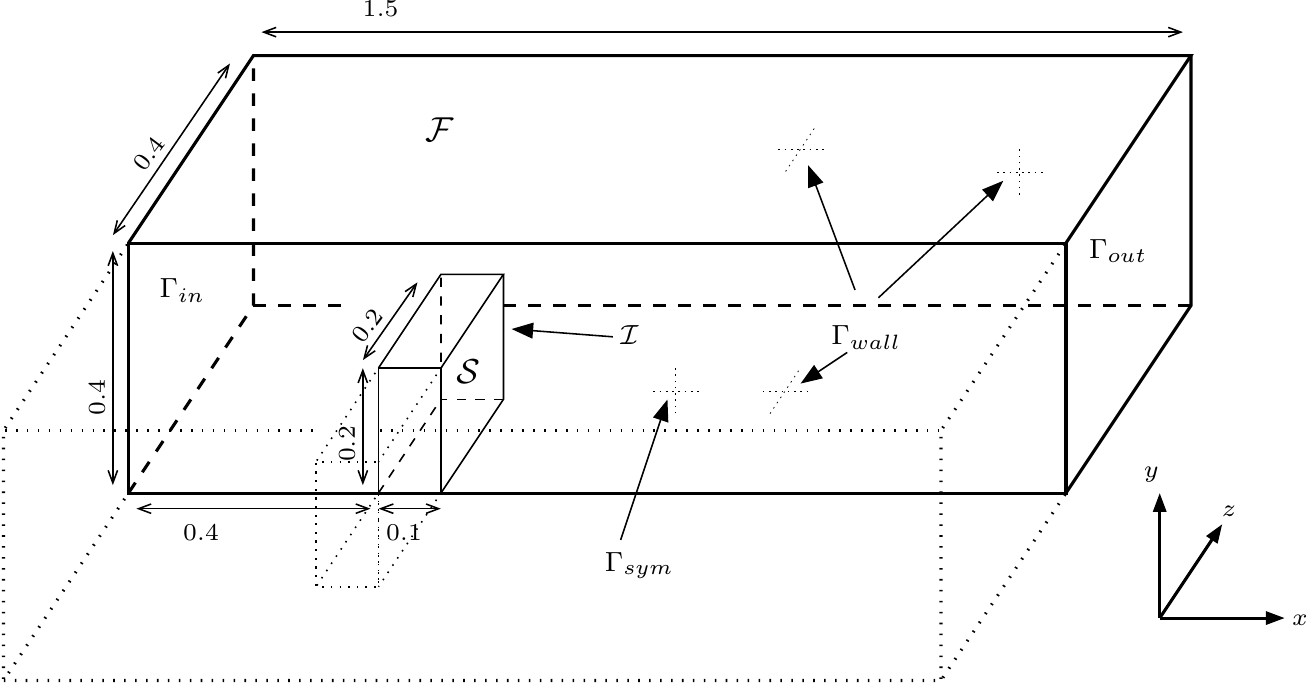}
    \end{center}
    \subcaption{3d configuration. Like in 2d with additional symmetry
      boundary $\Gamma_{sym}$, where $\vt\cdot\nt=0$ is prescribed. }
  \end{subfigure}
  \caption{Configuration of the test cases in 2d (top) and 3d
    (bottom). Besides the geometry in reference configuration we also
    specify the
    boundary conditions.}
  \label{fig:geometry}
\end{figure}

Figure~\ref{fig:geometry} shows the geometry of the test cases in the
2d and 3d configuration. Both test cases
assemble the flow around a wall-mounted elastic obstacle that will
undergo a deformation. The problem is driven by a parabolic
(bi-parabolic in 3d) inflow profile on $\Gamma_{in}$ which is
oscillating in time
\[
\begin{aligned}
  \vt_{in}^{2d}(t,\,y)&=s(t)\frac{y(H-y)}{(H/2)^2}\frac{3}{2}\bar\vt_{2d},\\
  \vt_{in}^{3d}(t,\,y,\,z)&=s(t)\frac{y(H-y)(H^2-z^2)}{(H/2)^2H^2}\frac{9}{8}\bar\vt_{3d},\qquad s(t)=\frac{1}{2}\big(1-\cos(\pi t)\big).
\end{aligned}
\]
By $H_{3d}=0.4$ and $H_{2d} = 1$ we denote the diameter of the flow domain and by
$\bar\vt_{2d}$ and $\bar\vt_{3d}$ the average flow rates. The scalar
function $s(t)$ oscillates with the period $T=\unit[1]{s}$.
In both cases the time interval is $I=[0,\,8]$.
All further parameters of this benchmark problem as shown in
Table~\ref{tab:parameters}. The configurations yield the maximum
Reynolds number $Re_{2d}=30$ in 2d and $Re_{3d}=20$ in 3d,
where $L_{2d}=0{.}5$ and $L_{3d}=0{.}2$, the height of the obstacles
is chosen as characteristic length scale. On the wall boundary
$\Gamma_{wall}$ we prescribe homogeneous Dirichlet conditions $\vt=0$
and on the outflow boundary $\Gamma_{out}$ the do-nothing condition,
see~(\ref{variational}). In 3d, we split the domain and introduce a
symmetry boundary, where we prescribe a free-slip, no-penetration
condition $\vt\cdot\nt=0$ and $\sigmat\nt\cdot\vec t=0$, where $\vec t$
are the tangential vectors.

\begin{table}[t!]
  \begin{center}
    \begin{tabular}{lll}
      \toprule
      Problem configuration    & 2d & 3d \\
      \midrule
      Fluid density            & $10^3\,\mathrm{kg}\cdot \mathrm{m}^{-3}$
      & $10^3\,\mathrm{kg}\cdot \mathrm{m}^{-3}$    \\
      Kinematic viscosity      & $2\cdot 10^{-2}\,\mathrm{m}^2 \mathrm{s}^{-1}$
      & $10^{-2}\,\mathrm{m}^2 \mathrm{s}^{-1}$     \\
      Average inflow velocity  & $1.2\,\mathrm{m}\cdot\mathrm{s}^{-1}$
      & $1.0\,\mathrm{m}\cdot\mathrm{s}^{-1}$       \\
      Solid density            & $10^3\,\mathrm{kg}\cdot \mathrm{m}^{-3}$
      & $10^3\,\mathrm{kg}\cdot \mathrm{m}^{-3}$    \\
      Shear modulus            & $1\cdot 10^6\,\mathrm{kg}\cdot \mathrm{m}^{-1}\cdot \mathrm{s}^{-2}$
      & $5\cdot 10^5\,\mathrm{kg}\cdot \mathrm{m}^{-1}\cdot \mathrm{s}^{-2}$\\
      Poisson ratio & $\nu=0.4$ & $\nu=0.4$ \\
      \bottomrule
    \end{tabular}
  \end{center}
  \caption{Configuration of the test problems in two and three
    dimensions. We indicate the parameters defining the different
    cases. }
  \label{tab:parameters}
\end{table}

\subsection{2d example}

\begin{figure}[t]
  \parbox{0.49\linewidth}{\centering
    \includegraphics[width=\linewidth]{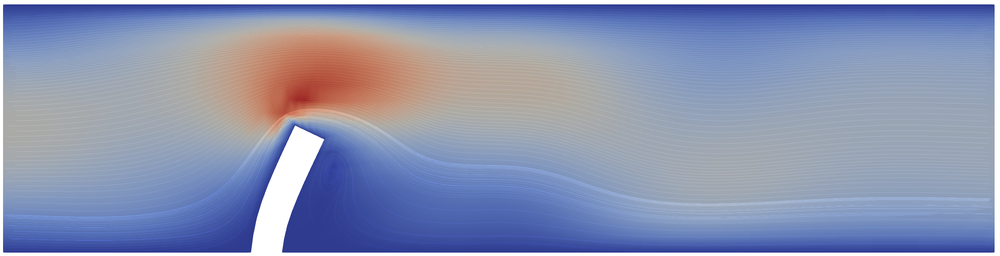}}
  \parbox{0.49\linewidth}{\centering
    \includegraphics[width=\linewidth]{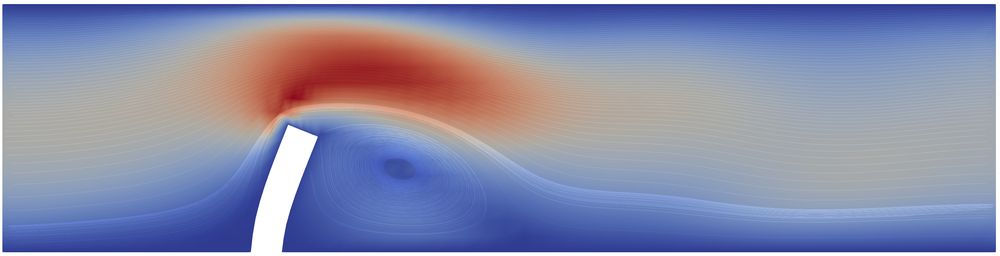}}

  \parbox{0.49\linewidth}{\centering
    \includegraphics[width=\linewidth]{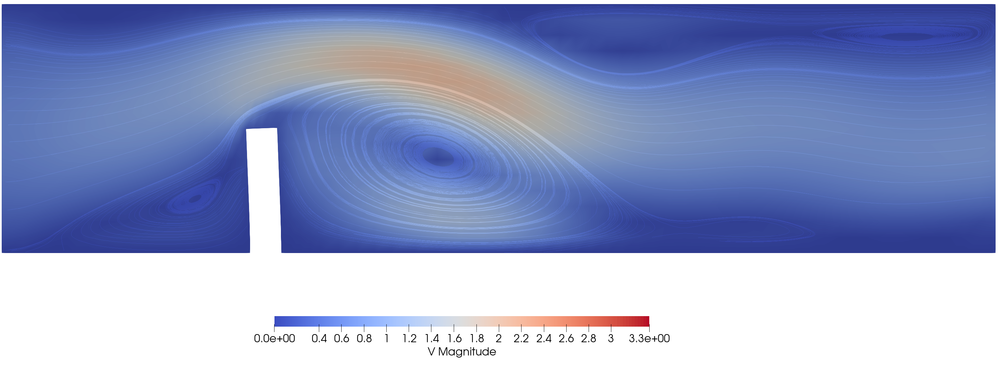}}
  \parbox{0.49\linewidth}{\centering
    \includegraphics[width=\linewidth]{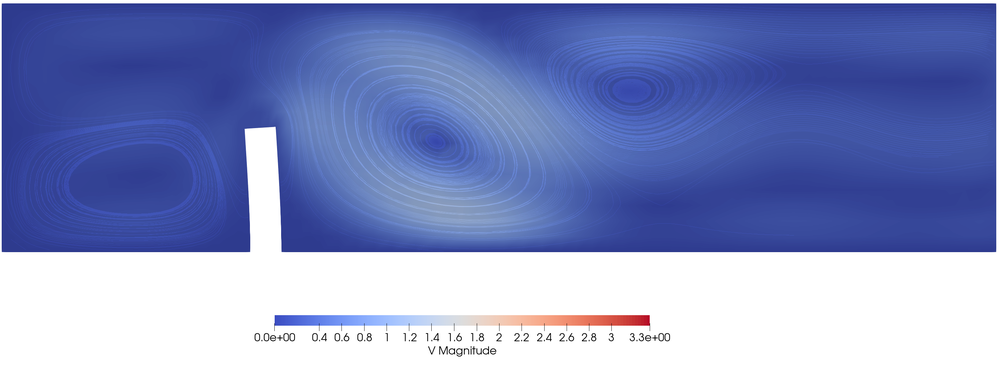}}
  \caption{Solutions on the last 4 subintervals (from left to right), the deformation is
    magnified by a factor of 2.}\label{fig:2dsolution}
\end{figure}
\edef\pgf@imagewidth{5cm}
\begin{figure}[t]
  \parbox{0.49\linewidth}{\centering \input{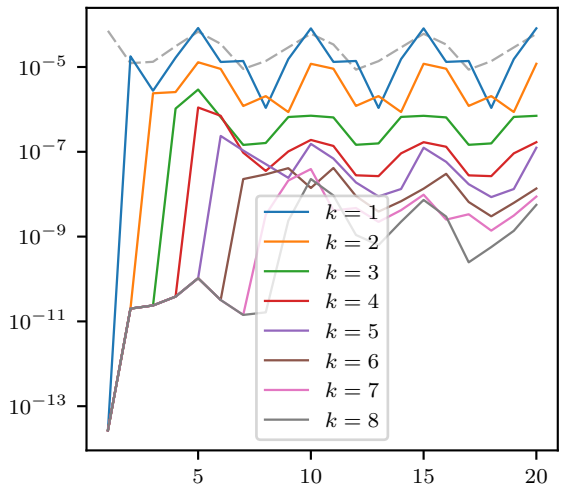}
    \subcaption{Coarse stepsize $K = 0.01$} }
  \parbox{0.49\linewidth}{\centering \input{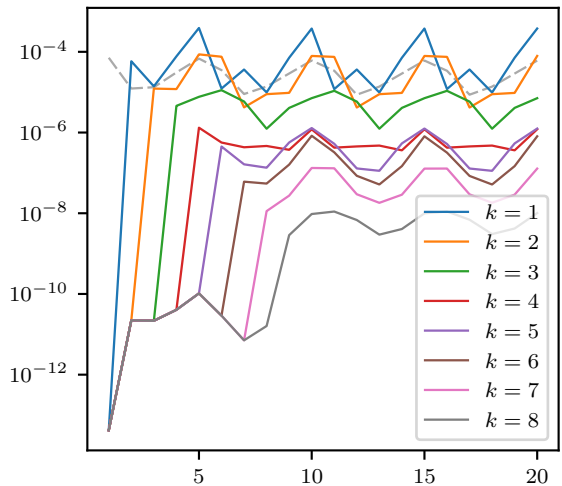}
    \subcaption{Coarse stepsize $K = 0.02$} }

  \parbox{0.49\linewidth}{\centering \input{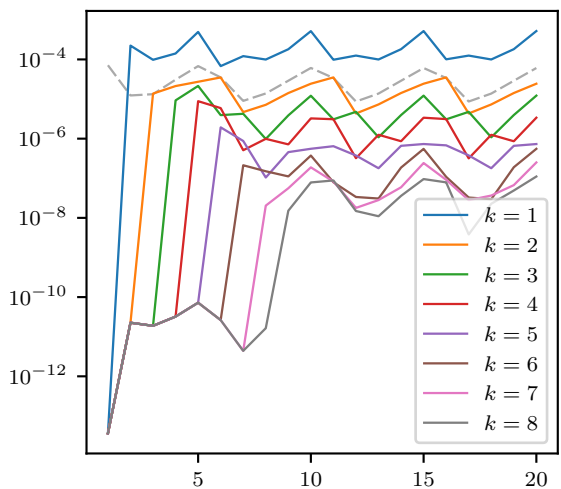}
    \subcaption{Coarse stepsize $K = 0.05$} }
  \parbox{0.49\linewidth}{\centering \input{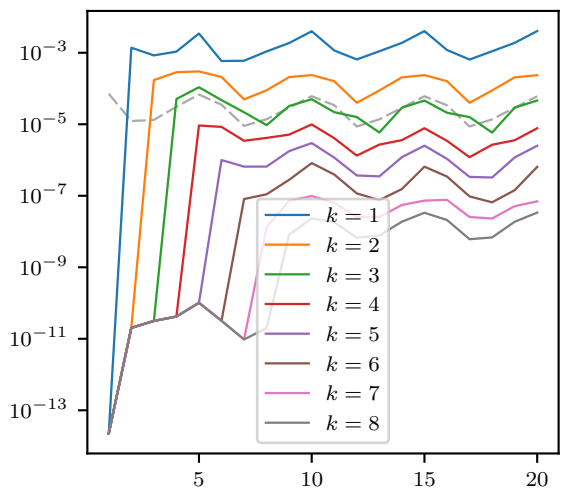}
    \subcaption{Coarse stepsize $K = 0.1$} }

  \caption{Velocity error between sequential (computed with the fine
    stepsize $k$) and Parareal solution with different coarse
    step sizes calculated at the interval boundaries as the relative
    error in the Frobenius norm
    $\frac{\lVert v_p-v_s\rVert_F}{\lVert v_s \rVert_F}$. Where $v_s$
    is the sequential solution and $v_p$ the Parareal solution. The
    dashed line is the difference between the sequential solution
    and a reference solution. All Parareal approximations below this
    line can be considered sufficiently accurate.}
  \label{fig:2derror}
\end{figure}
Figure~\ref{fig:2dsolution} shows the solutions of the problems at
time $t_{17}=6{.}8,\,t_{18}=7{.}2,\,t_{19}=7{.}6,\,t_{20}=8$. The dynamics of this problems is
enforced by the oscillating right hand side that causes a periodic
motion of the elastic obstacle.

In figure~\ref{fig:2derror} we show the convergence behavior of the
Parareal method for four choices of the large time step $K$. For each
iteration $k$ of the Parareal method we plot the relative error
between the Parareal approximation and a sequential solution (using
the same small step size).
The time interval is split into 20 subintervals and we clearly
recognize that after the \(k\)-th iteration, the solution on the
\(k\)-th interval does not get any better. All plots also include a dashed
line (this is exactly the same line in every plot) which indicates the
discretization error, e.g. the error between the sequential solution
with stepsize $k$ and a refined reference solution. If the Parareal
error is below this line, the Parareal approximation can be considered
sufficiently accurate, e.g. considering $K=0.01$, two iterations $k=2$
already yield a Parareal approximation error that is smaller than the
discretization error.

Furthermore, in figure~\ref{fig:speedup} we give an overview of the
speedup that is obtained in comparison to the sequential simulation.
In this example the optimal choice would be $K=0{.}05$, stopping after
3 iterations such that the discretization error is dominant and
resulting in a speedup of $3{.}45$.

That is below the theoretical value of $5{.}78$ estimated by
from~(\ref{eq:speedup}). The sub optimal performance is mostly due to the
strong nonlinearity of the fluid-structure interaction problem. The
nonlinear problems are approximated with a Newton's method that
benefits form small time-steps as they limit the role of the ALE
map. Using coarser time steps for the prediction (that should
theoretically improve the speedup) we require more Newton iterations
such that the efficiency is reduced. The theoretical predictions can
only be reached when
the computational costs of each time step are exactly the same,
which we cannot guarantee.

\FloatBarrier
\clearpage

\subsection{fsi-3 benchmark problem. A case where the algorithm fails}\label{subsec:fail}

\begin{figure}[h!]
  {\centering \input{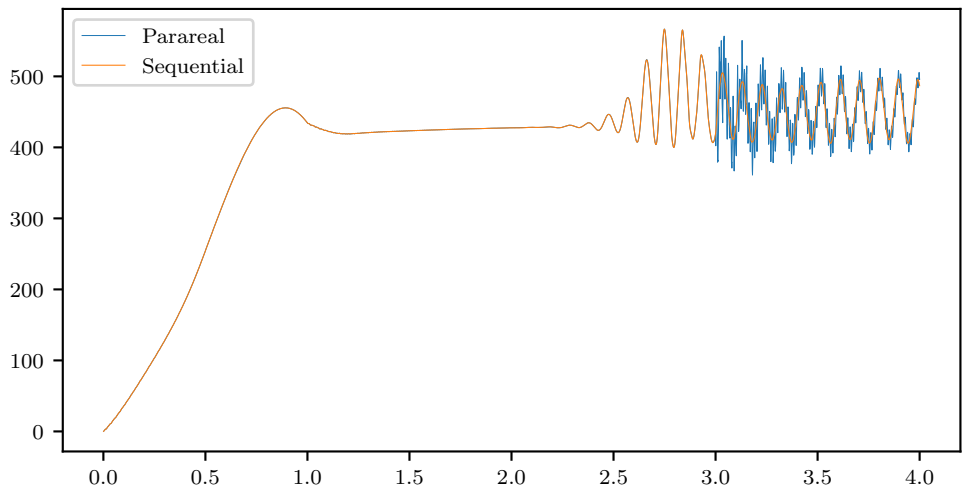} }
  \caption{fsi-3 benchmark problem (\cite{HronTurek2006}).
    Drag coefficient from the Parareal method after 2
    iterations in comparison to the sequential solution at
    time step size $k=10^{-3}$. The coarse stepsize was chosen as
    $k=10^{-2}$. Oscillations occur after each Parareal update.}\label{fig:fsi3}
\end{figure}

In the first experiments, we considered the {\itshape fsi-3 benchmark}
problem as introduced by Hron and Turek~\cite{HronTurek2006}.  In this
case Parareal did not converge. Instead we observed
pressure oscillations as shown in figure~\ref{fig:fsi3}. One possible
explanation for these instabilities that are also known from the
simulation of incompressible flow problems on moving
meshes~\cite{BesierWollner2012} is
the violation of the weak divergence freeness after the Parareal
update
$(\ut,\,\vt)_{m+k}^{\textrm{final}}=(\ut,\,\vt)_{m+k-1}^\PC +
(\ut,\,\vt)_{m+k}^\pf - (\ut,\,\vt)_{m+k}^\PC$, as the nonlinear ALE
formulation $(J\Ft^{-1}:\nabla\vt^T,\xi)_{\FL}$ with
separate updates in velocity and deformation does not conserve weak
incompressibility and gives rise to pressure oscillations
\begin{multline*}
  i=1,2:\quad \big( \operatorname{det}(I+\nabla
  \ut_i)(I+\nabla\ut_i)^{-1}:\nabla\vt^T_i,\xi\big)_{\FL}=0\\
  \not\Rightarrow\quad
  \big( \operatorname{det}(I+\nabla
  (\ut_1+\ut_2))
  (I+\nabla(\ut_1+\ut_2))^{-1}:\nabla(\vt_1+\vt_2)^T,\xi\big)_{\FL}=0.
\end{multline*}
Our approaches to avoid this problem,
e.\,g.\ applying a Stokes projection to the initial solutions on the
subintervals in every time steps, see~\cite{BesierWollner2012}, did
not give remedy.

Another and more likely reason for
the cause of these
instabilities is the dependency of the damping properties of the
shifted Crank-Nicolson time-stepping scheme on the
time-step. Depending on the time step size, the problem shows a
slightly shifted transient phase, e.g. the
dominant oscillations in the drag coefficient are slightly and when
the Parareal algorithm was applied these shifts caused instabilities.
If the coarse time step was chosen so small that the dynamics did not differ
in a destructive way, no speedup was achieved anymore. It could be
worthwhile to test a damping $\theta = \frac{1}{2}+\alpha(k) k$ with a
special choice of the parameter $\alpha(k)$, such that large time step and
short time step have the same dissipation properties.
In~\cite{ruprecht2018wave} Ruprecht suggests modifying the update in the parareal
algorithm, such that information about the dissipation properties is
taken into account.

\FloatBarrier

\subsection{3d example}

\begin{figure}[t]
  \parbox{0.49\linewidth}{\centering
    \includegraphics[width=\linewidth]{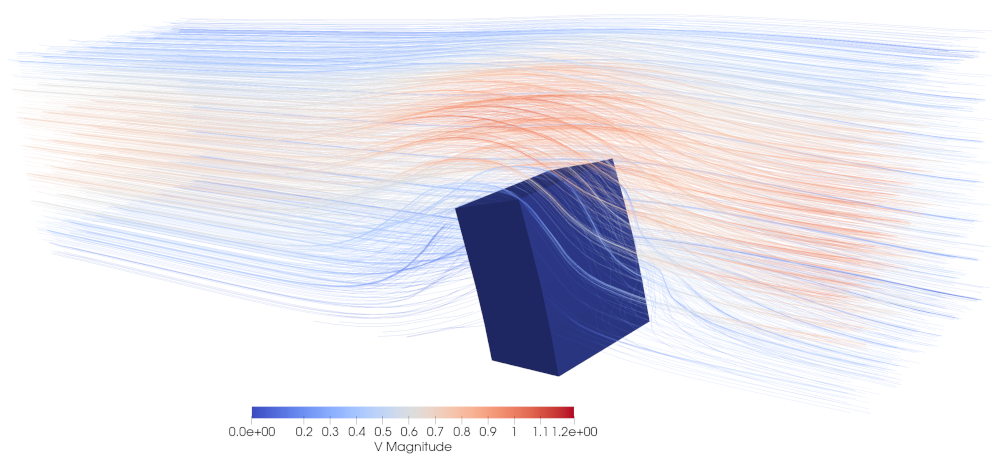}}
  \parbox{0.49\linewidth}{\centering
    \includegraphics[width=\linewidth]{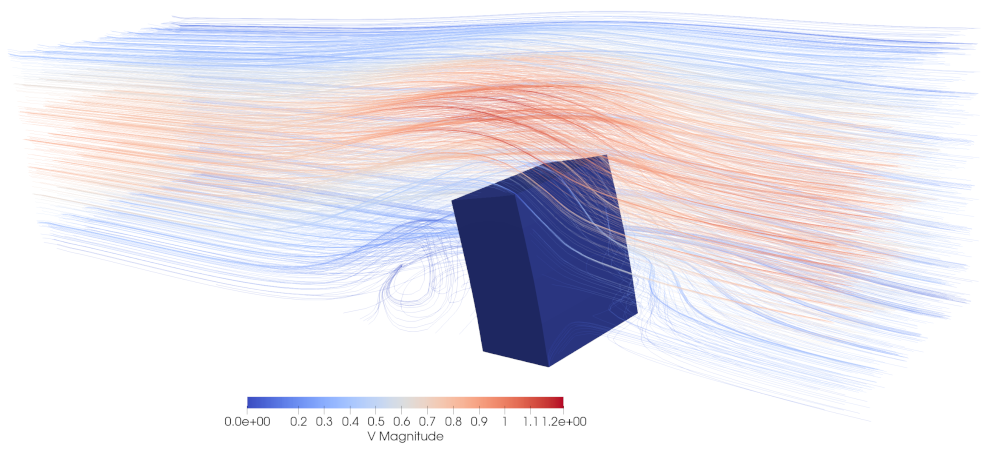}}

  \parbox{0.49\linewidth}{\centering
    \includegraphics[width=\linewidth]{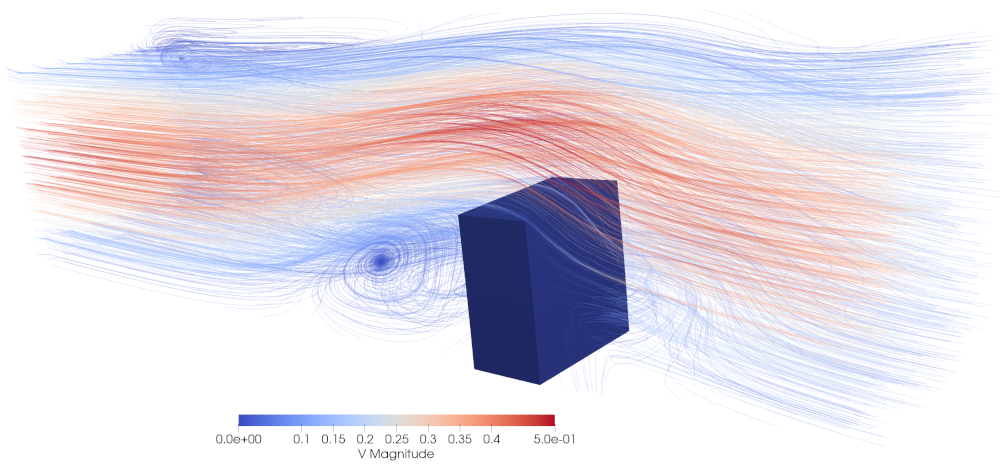}}
  \parbox{0.49\linewidth}{\centering
    \includegraphics[width=\linewidth]{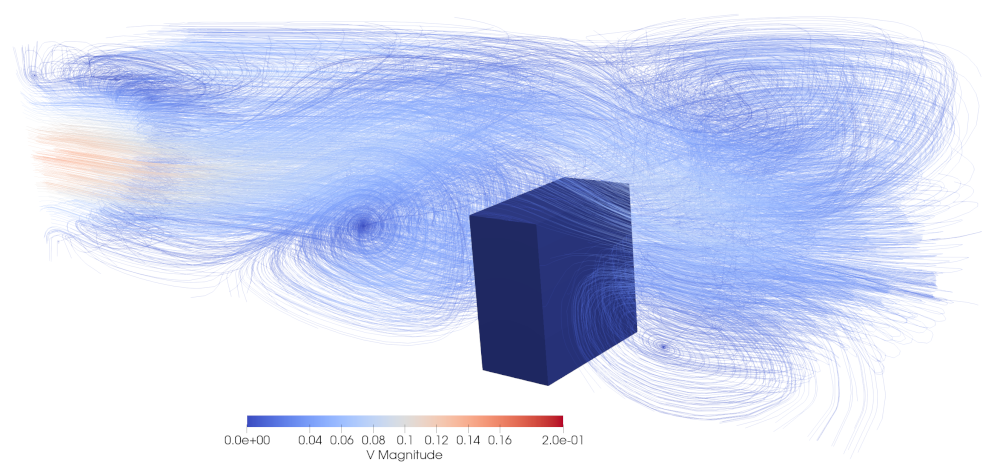}}
  \caption{Solutions on the last 4 subintervals (from left to right), the deformation is
    magnified by a factor of 10.}\label{fig:3dsolution}
\end{figure}

\edef\pgf@imagewidth{5cm}
\begin{figure}[t]
  \parbox{0.49\linewidth}{\centering \input{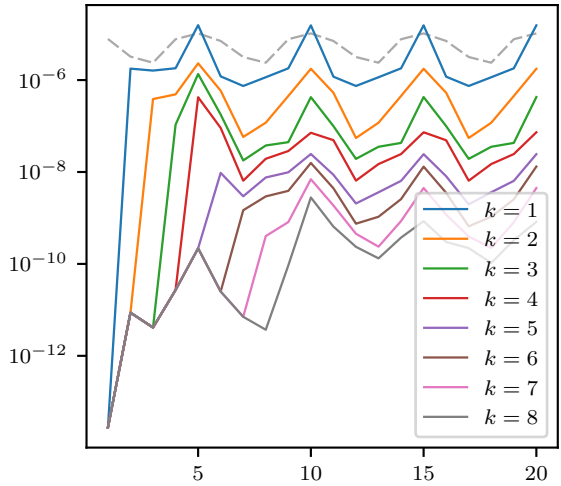}
    \subcaption{Coarse stepsize $K = 0.01$} }
  \parbox{0.49\linewidth}{\centering \input{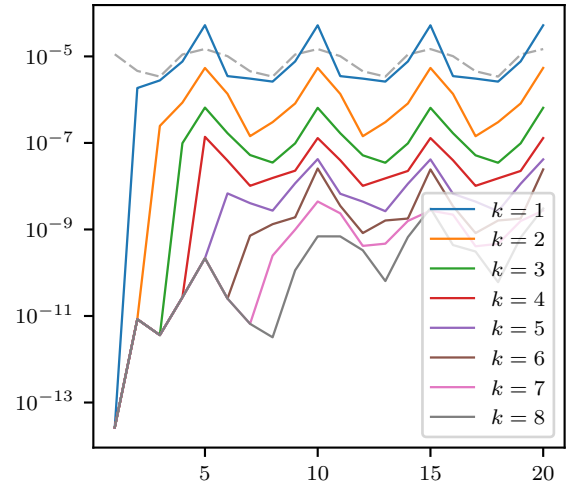}
    \subcaption{Coarse stepsize $K = 0.02$} }

  \parbox{0.49\linewidth}{\centering \input{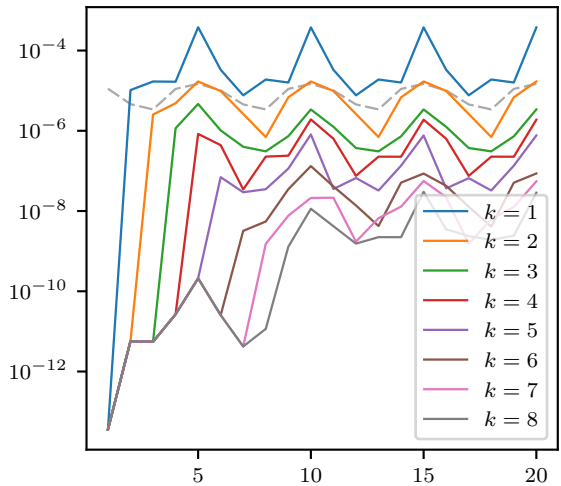}
    \subcaption{Coarse stepsize $K = 0.05$} }
  \parbox{0.49\linewidth}{\centering \input{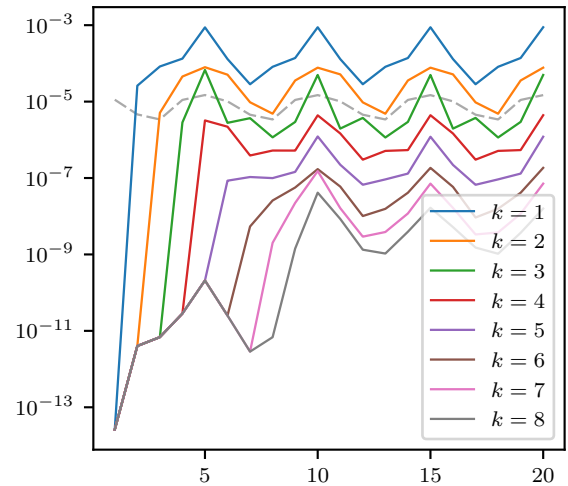}
    \subcaption{Coarse stepsize $K = 0.1$} }

  \caption{Velocity error from the 3d case. We see the error between
    sequential (computed with the fine
    stepsize $k$) and Parareal solution with different coarse
    step sizes. They are calculated at the subinterval boundaries as the relative
    error in the Frobenius norm
    $\frac{\lVert v_p-v_s\rVert_F}{\lVert v_s \rVert_F}$. Where $v_s$
    is the sequential solution and $v_p$ the Parareal solution. The
    bright grey line is the difference between the sequential solution
    and a reference solution. }\label{fig:3derror}
\end{figure}

\begin{figure}[t]
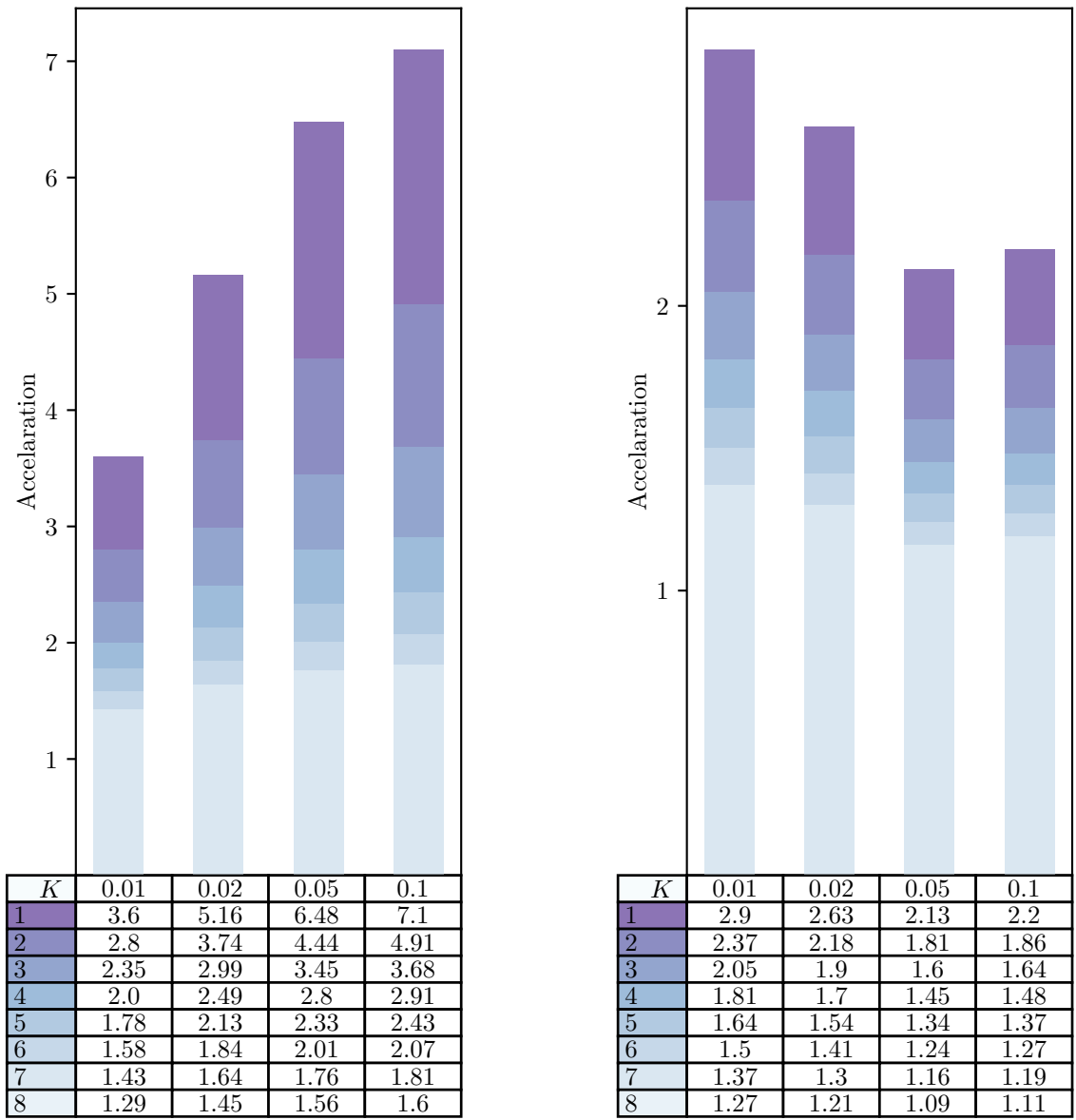

  \parbox{0.49\linewidth}{\centering
    \input{Plots/wm_2d/timelog20_2d.pgf} \subcaption{The 2d case} }
  \parbox{0.49\linewidth}{\centering
    \input{Plots/box3d/timelog20_3d.pgf} \subcaption{The 3d case} }
  \caption{Speedup of Parareal using 20 cores.}\label{fig:speedup}
\end{figure}

Figure~\ref{fig:3dsolution} shows the solution of the 3d problem
$t_{17}=6{.}8,\,t_{18}=7{.}2,\,t_{19}=7{.}6,\,t_{20}=8$, analogous to the 2d case. In
figure~\ref{fig:3derror} we show the relative errors of the Parareal
approximation in comparison to a sequential simulation for each
iteration. Again, a dashed line represents the discretization
error of the problem by comparison to a refined reference solution. We
observe a similar behavior like in the 2d case: The Parareal error is
quickly below the level of the discretization error and initial
intervals do not benefit from further steps. However,
figure~\ref{fig:speedup} showing the speedups indicates a
substantially reduced performance of the Parareal scheme when applied
to this 3d test case.
The main cause for this deterioration is the linearization of the
problem: Newton's method takes longer to converge for
big timesteps and the corrections in the Parareal algorithm only
seem to aggravate this problem.
Although the instabilities noted in~\ref{subsec:fail} could not be
observed here, they start to rise: Larger timesteps influence the
dynamics and if they get too large, the problems are not feasible in a
stable manner anymore.

Due to these problems the Parareal algorithm is not guaranteed
to give beneficial speedups in the 3d case. By choosing the timesteps
carefully a speedup of $2{.}18$ can still be obtained.

\section{Conclusion}
We applied the Parareal algorithm to fluid-structure interaction
problems in 2d and 3d. Two cases were analyzed, in which speedup and
convergence in very few iterations were achieved.
We also saw notable limitations of the algorithm when applied to
fluid-structure-interactions. One reason is the hyperbolic character
of the solid problem that prevents the use of too coarse
timesteps, as they can result in a different dynamics and an offset
between predictor and corrector. The fsi-3 benchmark problem by
Hron and Turek showed to be too challenging and it was not possible to
pick discretization parameters that gave speedup and a stable solution
at the same time.

Another issue is the severe
nonlinearity of fluid-structure
interactions. The design on nonlinear (and linear) solvers is already
a challenge and coarse timesteps lead to an increased effort that limits
the possible speedup.

The Parareal algorithm needs careful investigation of the problem
before it can be successfully applied. Choosing the right time step sizes
is crucial and even then it could still fail.
If the problem is suitable on the other hand, it offers
a simple way to speedup the solution.

\section*{Acknowledgement}

The authors acknowledge the financial support by the Federal Ministry of
Education and Research of Germany, grant number 05M16NMA as well as
the GRK 2297 MathCoRe, funded by the Deutsche
Forschungsgemeinschaft, grant number 314838170.

\end{document}